\documentclass[11pt,a4paper,english,reqno]{amsart}
\usepackage{amsmath,amssymb,amsfonts,epsfig,mathrsfs}
\usepackage[T1]{fontenc}

\usepackage{color}
\usepackage{array}
\usepackage{amsthm}
\usepackage{amstext}
\usepackage{graphicx}
\usepackage{setspace}
\usepackage[margin=2.5cm]{geometry}
\usepackage{bbm}
\usepackage{color}
\usepackage{enumitem}
\usepackage{undertilde}
\allowdisplaybreaks[4]

\usepackage{pgfplots}

\usepackage{amscd,psfrag}
\usepackage{yhmath}
\usepackage[mathscr]{eucal}

\usepackage{nicefrac}

\usepackage{slashed}

\makeatletter
\pdfpageheight\paperheight
\pdfpagewidth\paperwidth

\usepackage{mathrsfs}

\usepackage{epstopdf}
\usepackage{chngcntr}
\counterwithin{figure}{section}
\usepackage{mathrsfs}

\usepackage{indentfirst}	

\usepackage[normalem]{ulem}
\theoremstyle{plain}

\newtheorem{definition}{Definition}[section]
\newtheorem{theorem}[definition]{Theorem}
\newtheorem*{theorem*}{Theorem}

\newtheorem*{remark*}{Remark}
\newtheorem*{sideremark*}{Side Remark}\newtheorem*{mt*}{Main Theorem}

\newtheorem*{claim*}{Claim}
\newtheorem*{q*}{Question}

\newtheorem*{corollary*}{Corollary}
\newtheorem*{proposition*}{Proposition}

\newcommand{\R}{\mathbb{R}}

\newcommand{\na}{\nabla}

\newcommand{\dd}{{\rm d}}

\newcommand{\e}{\epsilon}

\newcommand{\map}{\rightarrow}

\newcommand{\M}{\mathcal{M}}

\newcommand{\can}{{\mathbf{can}}}

\allowdisplaybreaks[4]

\def\XXint#1#2#3{{\setbox0=\hbox{$#1{#2#3}{\int}$ }
\vcenter{\hbox{$#2#3$ }}\kern-.6\wd0}}

\numberwithin{equation}{section}
\numberwithin{figure}{section}

\title{Counterexamples to the $L^p$-Calder\'{o}n--Zygmund Estimate on Open Manifolds}

\author{Siran Li}
\address{Siran Li: Department of Mathematics, Rice University, MS 136
P.O. Box 1892, Houston, Texas, 77251-1892, USA\, $\bullet$ \,  Department of Mathematics, McGill University, Burnside Hall, 805 Sherbrooke Street West, Montreal, Quebec, H3A 0B9, Canada.}

\email{\texttt{Siran.Li@rice.edu}}

\keywords{Calder\'{o}n--Zygmund Estimate; Counterexample; Open Manifold; Warped Product; Laplace--Beltrami; Hessian.}
\subjclass[2010]{Primary: 	58J05, 53C21; Secondary: 35J05, 35J08.}

\date{\today}

\begin{document}

\maketitle

\begin{abstract}
Based on a construction due to B.  G\"{u}neysu and S. Pigola (\textit{Adv. Math.} \textbf{281} (2015), pp.353--393), for each $p \in [1,\infty]$ and $m \in \mathbb{Z}_{\geq 2}$, we exhibit an $m$-dimensional Riemannian open manifold $\M$ on which the $L^p$-Calder\'{o}n--Zygmund estimate
\begin{equation*}
\|\na \na f\|^p_{L^p} \leq C_1 \|\Delta f\|_{L^p}^p + C_2 \| f\|_{L^p}^p \qquad \text{ for all } f \in C^\infty_c(\M) 
\end{equation*}
is false for any $C_1, C_2$ depending on $m$ and $p$. Therefore, one must impose further geometric conditions on the manifold to ensure the validity of the Calder\'{o}n--Zygmund estimate.
\end{abstract}

\section{Introduction}
In this short note, we prove the following:
\begin{theorem}\label{thm}
For any $1\leq p \leq \infty$, any $m \in \mathbb{N}$ and any positive constants $C_1,C_2$, there is an open Riemannian manifold $(\M,g)$ of dimension $m$ such that the $L^p$-Calder\'{o}n--Zygmund estimate is invalid. More precisely, there is a smooth function $f:\M \map \R$ such that
\begin{equation}\label{E}
\|\na \na f\|^p_{L^p} > C_1 \|\Delta f\|_{L^p}^p + C_2 \| f\|_{L^p}^p.
\end{equation}
\end{theorem}

Throughout, a Riemannian manifold (without boundary) $(\M,g)$ is said to be open if it is non-compact and geodesically complete. The $2$-tensor field $\na_g \na_g f$ denotes the Hessian of $f$. Its trace under the metric is the Laplace--Beltrami operator, denoted by $\Delta_g f$. In this note, all the Hessian and Laplace--Beltrami operators shall be taken with respect to a fixed metric $g$ (the one in \eqref{metric} below); let us abbreviate by $\na \na f :=\na_g\na_g f$ and $\Delta f :=\Delta_g f$.

In the Euclidean space $\R^m$, the classical estimate 
\begin{equation}\label{cz}
\|\na \na f\|^p_{L^p} \leq C_1 \|\Delta f\|_{L^p}^p + C_2 \| f\|_{L^p}^p\qquad \text{ for all } f \in C^\infty_c(\R^m) \text{ and every } p \in ]1,\infty[
\end{equation}
was established by Calder\'{o}n--Zygmund in the seminal paper \cite{cz}. Here the constants $C_1$, $C_2$ depend only on $p$ and $m$. A natural question is the validity of \eqref{cz} on a Riemannian manifold $(\M,g)$. Many works are devoted to proving \cite{cz} on $(\M,g)$ which satisfies certain geometric assumptions, {\it e.g.}, the boundedness of Ricci or sectional curvatures, the boundedness of the injectivity radius away from zero, and the doubling property for the Riemannian volume measure. We refer to Cheeger--Gromov--Taylor \cite{cgt}, Strichartz \cite{s}, Taylor \cite{t}, Wang \cite{w} and G\"{u}neysu--Pigola \cite{gp} for details; also see the many references cited therein.

On the other hand, in \cite{gp} G\"{u}neysu--Pigola constructed a 2-dimensional complete manifold $(\M,g)$ on which \eqref{cz} is invalid $p=2$. To the author's knowledge, this is among the first negative results for the Calder\'{o}n--Zygmund estimates.  Our goal here is to generalise the arguments in \cite{gp} to prove Theorem \ref{thm} for the whole range of indices $p \in [1,\infty]$ and $m \in \mathbb{Z}_{\geq 2}$.

Before starting the proof, let us make three remarks:
\begin{enumerate}
\item
Throughout this paper, $\|\bullet\|_{L^p}$ denotes the $L^p$-norm of tensor fields on $\M$, taken with respect to the metric $g$ in \eqref{metric} below. For the definition and discussions on Sobolev spaces over manifolds, see {\it e.g.} Hebey \cite{h}. 
\item
The Calder\'{o}n--Zygmund estimate is known to be false for $p=1$ and $p=\infty$ on Euclidean spaces; see Ornstein \cite{o} and McMullen \cite{m}. It thus remains to prove for $1<p<\infty$.
\item
Our proof is crucially based on the construction in Theorem B, \cite{gp} by G\"{u}neysu--Pigola.  
\end{enumerate}

\section{Construction of the Manifold $(\M,g)$}

In this section, we construct the manifold $(\M,g)$ which leads eventually to the proof of Theorem \ref{thm}. The presentation in this section works for all $m \in \mathbb{Z}_{\geq 2}$. It involves the choice of several parameters, which will be specified in subsequent sections. 

\smallskip
\noindent
{\bf Warped Product.} The manifold $\M$ we choose is the Euclidean space $\R^m$ equipped with the warped product manifold:
\begin{equation}\label{metric}
g= dr\otimes dr + \sigma^2(r)\, \can^{m-1},
\end{equation}
where $\can^{m-1}$ is the canonical round metric $(m-1)$-dimensional Euclidean sphere.
It is known as a {\em warped product} manifold. Note that the space forms are special examples of warped products: $\M=\mathbb{S}^m$ if $\sigma=\sin$, $\M=\R^m$ if $\sigma={\bf Id}$, and $\M = \mathbb{H}^m$ if $\sigma = \sinh$. We shall choose the {\em warping function} $\sigma$ to be non-negative, smooth and growing to infinity as the radial coordinate $r \nearrow +\infty$. Thus $(\M, g)$ is an open manifold. The warped products are central objects of many recent works on geometric analysis; {\it cf. e.g.} \cite{gl} by Guan--Lu.

\smallskip
\noindent
{\bf Green's Function.} Let $\widetilde{G}(x)$ be the Green's function of the Laplace--Beltrami operator on $\M$ as above. Since $g$ in \eqref{metric} is rotationally symmetric, there is a function $G:[0,\infty[\map \R$ such that $\widetilde{G}(x) = G(r)$ for $r:=|x|$. Writing $\Delta$ in polar coordinates, we find that
\begin{equation}\label{green fn}
\Delta \widetilde{G} = 0 \quad \Longleftrightarrow \quad G'' + (m-1) \frac{\sigma'}{\sigma}G = 0\qquad \text{ on } \M \sim \{0\}.
\end{equation}

\smallskip
\noindent
{\bf Hessian and Laplacian.} The key idea of the construction, as in \cite{gp}, is to take $f$ to be a localised version of the Green's function $G$. For $k \in \mathbb{N}$ and $[\alpha_k, \beta_k] \subset \R$, let $\phi_k \in C^\infty_c ([\alpha_k, \beta_k])$ be a cut-off function.  Here $\phi_k, \alpha_k, \beta_k$ will be specified later. Then, define
\begin{equation}\label{uk}
u_k(r) := \phi_k \circ G(r).
\end{equation}
In the end, one shall set $f:=u_k$ for some sufficiently large $k$. 

Direct computations in \cite{gp} lead to the following formulae for the Hessian and the Laplace--Beltrami of $u_k$, as well as the volume density of $g$:
\begin{eqnarray}
&& \na \na u_k = u_k''\,dr\otimes dr + \sigma {\sigma'} u_k'\, \can^{m-1}, \label{hessian}\\
&& \Delta u_k = u_k'' + (m-1)\frac{\sigma' u_k'}{\sigma},\label{laplacian}\\
&& \sqrt{\det \, g} = \sigma^{m-1}.
\end{eqnarray}
Throughout $\sigma, G, u_k$ are functions of $r$ only; $\sigma'$, $u_k''$ etc. denote the derivatives in $r$.

In the rest of this section, fixing a $p\in ]1,\infty[$, we collect some preliminary estimates for the $L^p$-norm of $u_k$, $\Delta u_k$ and $\na \na u_k$. First of all, neglecting the radial components in \eqref{hessian}, we have 
\begin{equation*}
|\na \na u_k|^p \geq \Big|u_k' \frac{\sigma'}{\sigma}\Big|^p.
\end{equation*}
Hence, denoting by $\gamma_m := {\rm Vol}_{\can^{m-1}}(\mathbb{S}^{m-1})$, the area of the unit sphere, we deduce
\begin{align*}
\|\na \na u_k(r)\|_{L^p} &= \gamma_m \bigg\{ \int_0^\infty \big|\na \na u_k\big|^p \sigma^{m-1}(r)\,\dd r \bigg\}^{\frac{1}{p}}\nonumber\\
&\geq \gamma_m \bigg\{ \int_0^\infty \Big|\phi_k'\big(G(r)\big) \, G'(r) \,\Big(\frac{\sigma'}{\sigma}\Big)(r)\Big|^p \sigma^{m-1}(r)\,\dd r \bigg\}^{\frac{1}{p}}\nonumber\\
&= \gamma_m \bigg\{ \int_0^\infty \Big|\phi_k'\big(G(r)\big)\Big|^p \Big|\frac{\sigma'}{\sigma}(r)\Big|^p  \sigma^{(1-m)(p-2)}(r) G'(r)\,\dd r \bigg\}^{\frac{1}{p}}.
\end{align*}
In the last line we used the identity $$G'(r)=\sigma^{1-m}(r).$$ A change of variable $r \mapsto s=G(r)$ yields that
\begin{equation}\label{hessian, Lp}
\|\na \na u_k\|_{L^p} \geq \gamma_m \bigg\{ \int_{\alpha_k}^{\beta_k} |\phi_k'(s)|^p \, \Big| \frac{\sigma'}{\sigma}\circ G^{-1}(s)\Big|^p \big[\sigma \circ G^{-1}(s)\big]^{(1-m)(p-2)} \,\dd s\bigg\}^{\frac{1}{p}}.
\end{equation}
For the Laplace--Beltrami, it is crucial to observe that
\begin{equation}\label{laplacian, 2}
\Delta u_k(r) = \phi_k'' \circ G(r) \, \sigma^{2-2m}(r),
\end{equation}
thanks to the defining property \eqref{green fn} of the Green's function. Thus, we have
\begin{align}\label{laplacian, Lp}
\|\Delta u_k\|_{L^p} = \gamma_m \bigg\{ \int_{\alpha_k}^{\beta_k} |\phi_k''(s)|^p \big[\sigma\circ G^{-1}(s)\big]^{2(p-1)(1-m)} \,\dd s\bigg\}^{\frac{1}{p}} 
\end{align}
Finally, note that
\begin{equation}\label{uk, Lp}
\|u_k\|_{L^p} = \gamma_m \bigg\{ \int_{\alpha_k}^{\beta_k} |\phi_k(s)|^p \big[\sigma \circ G^{-1}(s)\big]^{2(m-1)}\,\dd s \bigg\}^{\frac{1}{p}}
\end{equation}

The key observation: Only the norm of $\sigma$ is involved in the upper bounds for $\|\Delta u_k\|_{L^p}$ and $\|u_k\|_{L^p}$, while $\sigma'$ is present in the lower bound for $\|\na \na u_k\|_{L^p}$; see \eqref{hessian, Lp}, \eqref{laplacian, Lp} and \eqref{uk, Lp}. As a consequence, by carefully choosing a highly oscillatory profile for $\sigma$, we may force $\|\na \na u_k\|_{L^p}$ to be much larger than $\|\Delta u_k\|_{L^p}$ and $\|u_k\|_{L^p}$, thus contradicting the Calder\'{o}n--Zygmund inequality.

\section{Proof for $m=2$}\label{sec: m=2}

In this section we prove Theorem \ref{thm} for $m=2$ by specifying the warping function $\sigma$. The proof is essentially an adaptation of the arguments for Theorem B in \cite{gp} by G\"{u}neysu--Pigola, which corresponds to the case $m=2$, $p=2$. For the sake of completeness, we shall explain in detail why the constructions in \cite{gp} works for all $p \in ]1,\infty[$.

First, we set $\alpha_k = k$ and $\beta_k = k+1$ for each $k \in \mathbb{N}$. 

Next, let us require the warping function $\sigma$ to satisfy the following:
\begin{equation}\label{sigma, first properties}
\begin{cases}
\sigma^{(2k)}(0)=0 \text{ for each } k \in \mathbb{N};\\
\sigma'(0)=1;\\
\sigma(t)>0 \text{ for any } t>0;\\
t\leq \sigma(t) \leq t+1 \text{ for any }t \geq 1.
\end{cases}
\end{equation}
When $m=2$, one has the simple comparison results (see p.377 in \cite{gp}):
\begin{equation}\label{comparison 1}
\log\Big(\frac{t+1}{2}\Big) \leq G(t) \leq \log t\qquad \text{ for all } t>1
\end{equation}
and
\begin{equation}\label{comparison 2}
e^s \leq \sigma \circ G^{-1} (s) \leq 2e^s\qquad \text{ for all } s>0.
\end{equation}
Moreover, there exists a {\em universal} constant $\delta>0$ such that for all sufficiently large $k$, we can find $h=h(k)> k$ such that 
\begin{equation}\label{h def}
[h,h+1] \subset \big[G^{-1}(k+\delta),\, G^{-1}(k+1-\delta)\big].
\end{equation}

In addition, we choose the cut-off function $\phi_k$ in \eqref{uk} as follows: Fix some $\phi \in C^\infty_c (]0,1[)$ such that $\phi \equiv {\bf Id}$ on $[\delta, 1-\delta]$ and $\phi \leq 1$, and then set $$\phi_k(t) := \phi(t-k)$$ for each $k\in\mathbb{N}$. Here $\delta>0$ is the same constant as in \eqref{h def}. We shall fix $\phi$ once and for all; in particular, $\|\phi\|_{C^2([0,1])}$ is bounded by a universal constant.

We can deduce from \eqref{laplacian, Lp}, \eqref{uk, Lp} and \eqref{comparison 2} the following bounds:
\begin{align}\label{uk Lp upper bd}
\|u_k\|_{L^p}^p \leq  2 (\gamma_2)^p e^{2k+2}
\end{align}
and 
\begin{align}\label{laplacian Lp upper bd}
\|\Delta u_k\|_{L^p}^p \leq \frac{(\gamma_2)^p}{2(p-1)4^{p-1}}\|\phi''\|_{L^\infty [0,1]}^p e^{-2(p-1)k}. 
\end{align}
So it remains bound $\|\na \na u_k\|_{L^p}^p$ from below.

For this purpose, we shall further specify $\sigma$. Consider the cube
\begin{equation*}
Q_k := [k,k+1] \times  [k,k+1]\qquad \text{ for each } k \in \mathbb{N};
\end{equation*}
from the previous constructions, the graph of $\sigma$ is contained in $\bigcup_{k=0}^\infty Q_k$ (in fact, in the union of the upper-left corners of $Q_k$). For certain sequence $\{n_k\} \subset \mathbb{N}$ increasing to $+\infty$ as $k$ grows, we take 
\begin{equation*}
\e_k := \frac{1}{2n_k}.
\end{equation*}
Define $\mathfrak{S}_k$ on $[k,k+1]$ by the ``sawtooth'' function on p.378, \cite{gp}: 
\begin{equation*}
\mathfrak{S}_k(t):=
\begin{cases}
&k+2j\e_k + \frac{\e_k+1}{\e_k}(t- k -2j\e_k)\\
&\qquad\qquad\qquad \text{ on } \big[k+2j\e_k, k+(2j+1)\e_k\big]\,\,\text{ for each } j \in \{0,1,\ldots, n_k\},\\
&k+(2j+1)\e_k +1 + \frac{1-\e_k}{\e_k} \big(k+(2j+1)\e_k -t\big)\\
&\qquad\qquad\qquad \text{ on } \big[k+(2j+1)\e_k, k+2(j+1)\e_k\big]\,\,\text{ for each } j \in \{0,1,\ldots, n_k\}.
\end{cases}
\end{equation*}
Then, one defines $\sigma|[k,k+1]$ by smoothing out the corners of $\mathfrak{S}_k$. More precisely, for each $k \in \mathbb{N}$ we can take $\sigma \in C^\infty([k,k+1])$ such that $$\sigma = \mathfrak{S}_k \qquad \text{ on } [k,k+1] \sim \bigsqcup_{j=0}^{n_k} \Big[ k+2j\e_k - \e_k^{10},\, k+2j\e_k + \e_k^{10} \Big]$$ 
and that $\|\sigma\|_{C^3} \leq 2$ in each of the small intervals removed.

The idea for the construction of $\mathfrak{S}_k$ is clear: its graph (lying in the upper-left corner of the cube $Q_k$) is obtained by continuously concatenating $n_k$ copies of the following ``sawtooth unit'' with step-length $(2\e_k)$ --- in the first $\e_k$ it grows with constant gradient $\nicefrac{(\e_k+1)}{\e_k}$, and in the second $\e_k$ it decreases with constant gradient $\nicefrac{(1-\e_k)}{\e_k}$. In particular, in the second half of each sawtooth unit, the norm of the gradient is large.

With the above choice of $\sigma$, we can continue the lower bound \eqref{hessian, Lp} for the Hessian of $u_k$ as in below. First, by the definition of $\phi_k$ and  \eqref{comparison 2}, we have
\begin{align*}
\|\na\na u_k\|_{L^p}^p \geq (\gamma_2)^p 2^{-p}e^{-p(k+1)} \int_{k+\delta}^{k+1-\delta} \big|\sigma'\circ G^{-1}(s)\big|^p \,\big|\sigma\circ G^{-1}(s)\big|^{2-p}\,\dd s.
\end{align*}
Considering separately $p\geq 2$ and $p<2$ and using again \eqref{comparison 2}, one deduces
\begin{align*}
\|\na\na u_k\|_{L^p}^p \geq \min\big\{1,2^{2-p}\big\} (\gamma_2)^p 2^{-p}e^{-p(k+1)} \int_{k+\delta}^{k+1-\delta} \big|\sigma'\circ G^{-1}(s)\big|^p\,\dd s.
\end{align*}
For $m=2$ we have $G'(r)=\sigma^{-1}(r)$, hence 
\begin{equation*}
(G^{-1})'(s) = \frac{1}{G'[G^{-1}(s)]} = {\sigma[G^{-1}(s)]}.
\end{equation*}
It follows that 
\begin{align*}
\|\na\na u_k\|_{L^p}^p &\geq \min\big\{1,2^{2-p}\big\} (\gamma_2)^p 2^{-p}e^{-p(k+1)} \int_{k+\delta}^{k+1-\delta} \big|\sigma'\circ G^{-1}(s)\big|^p  \frac{1}{\sigma \circ G^{-1}(s)} (G^{-1})'(s) \,\dd s\\
&\geq \min\big\{1,2^{2-p}\big\} (\gamma_2)^p 2^{-p-1}e^{-p(k+1)}e^{-k-1+\delta} \int_{G^{-1}(k+\delta)}^{G^{-1}(k+1-\delta)} |\sigma'(r)|^p \,\dd r.
\end{align*}
Here we have used \eqref{comparison 1} once more.  

Recall that the universal constant $\delta$ is chosen right beneath \eqref{comparison 2}. For $k$ sufficiently large, we have selected $h=h(k)> k$ in \eqref{h def} so that
\begin{equation*}
\|\na\na u_k\|_{L^p}^p \geq \min\big\{1,2^{2-p}\big\} (\gamma_2)^p 2^{-p-1}e^{-(p+1)(k+1)+\delta} \int_{h}^{h+1} |\sigma'(r)|^p \,\dd r.
\end{equation*}
Thanks to the choice of $\sigma$, on $[h,h+1]$ the norm of the gradient $|\sigma'|$ is larger than $(2n_k-1)$ on more than $n_k$ intervals longer than $(\e_k - \e_k^{10})$, where $2n_k \e_k =1$. Thus,
\begin{align}\label{lower bound for hessian}
\|\na\na u_k\|_{L^p}^p &\geq \min\big\{1,2^{2-p}\big\} (\gamma_2)^p 4^{-p}e^{-(p+1)(k+1)+\delta} (2n_k-1)^p n_k (\e_k-\e_k^{10}) \nonumber\\
&\geq \min\big\{1,2^{2-p}\big\} (\gamma_2)^p 2^{-1-3p} e^{-(p+1)(k+1)+\delta} (1-\e_k^{9}) \e_k^{-p}.
\end{align}

We may now derive the contradiction by comparing \eqref{lower bound for hessian} with \eqref{uk Lp upper bd} and \eqref{laplacian Lp upper bd}. Note that
\begin{equation*}
\|u_k\|_{L^p}^p \lesssim e^{2k+2}\qquad \text{ and } \qquad \|\Delta u_k\|^p_{L^p} \lesssim e^{-2(p-1)k} \lesssim 1,
\end{equation*}
where the constants in $\lesssim$ depend on $p$, $C_1$, $C_2$ and $\|\phi''\|^p_{L^\infty([0,1])}$. 
On the other hand, 
\begin{equation*}
\|\na\na u_k\|_{L^p}^p \gtrsim e^{-(p+1)k} (1-\e_k^{9})\e_k^{-p}.
\end{equation*}
By further requiring for large $k\in\mathbb{N}$ that $\e_k \leq 100^{-1}$, we get $$\|\na\na u_k\|_{L^p}^p \gtrsim e^{-(p+1)k}\e_k^{-p},$$ with the constant in $\gtrsim$ depends only on $p$. Therefore, we can achieve \eqref{E} by choosing {\it e.g.},  $$\e_k := Ce^{-e^{k}}$$ for a suitable constant $C=C(p,C_1,C_2,\|\phi\|_{C^2([0,1])})$. Thus, choosing $k$ to be sufficiently large, we can complete the proof of Theorem \ref{thm} for $m=2$.

\section{Proof for $m \geq 3$}\label{sec: mD}

In this section, we prove Theorem \ref{thm} for arbitrary $m \geq 3$. 

The new feature is that the cubes $Q_k$ in $\S \ref{sec: m=2}$ are not available, since we cannot choose the warping function to satisfy $t\leq \sigma(t) \leq t+1$ for all $t \geq 1$. Instead, we shall choose an infinite sparse family of cubes $\{Q'_k\}$ sandwiched between the graphs of $t \mapsto t^{\frac{1}{m-1}}$ and $t \mapsto (t+1)^{\frac{1}{m-1}}$. Necessarily the size of the $Q'_k$ will shrink to zero as $k \nearrow \infty$; nevertheless, we can prescribe the rate of oscillation of $\sigma$ to be much larger than the shrinking rate of $Q'_k$. This is enough to conclude Theorem \ref{thm} for $m \geq 3$.

Now we start the proof. First of all, let us observe that the estimates \eqref{hessian, Lp}, \eqref{laplacian, Lp} and \eqref{uk, Lp} are valid for all $m \in \mathbb{Z}_{\geq 2}$, and that the radial Green's function again verifies
\begin{equation*}
G(r) = \int_1^r \sigma^{1-m}(t)\,\dd t.
\end{equation*}
We shall pick a $\sigma$ satisfying $G(+\infty)=+\infty$, which ensures the parabolicity of $(\M,g)$. For brevity we write
\begin{equation*}
\alpha \equiv \alpha_m := \frac{1}{m-1}.
\end{equation*}
Then, we choose $\sigma$ to satisfy a set of properties similar to those in \eqref{sigma, first properties}:\begin{equation}\label{sigma, md}
\begin{cases}
\sigma^{(2k)}(0)=0 \text{ for each } k \in \mathbb{N};\\
\sigma'(0)=1;\\
\sigma(t)>0 \text{ for any } t>0;\\
t^{\alpha}\leq \sigma(t) \leq (t+1)^\alpha \text{ for any }t \geq 1.
\end{cases}
\end{equation}
The motivation is to require the norm of $\sigma$ to be comparable to $t^\alpha$ without introducing a singularity at the origin. This can be achieved, {\it e.g.}, by gluing $\sigma|[1,\infty[$ to $\sinh$ or $\sin$ near $r=0$.

Notice that \eqref{comparison 1} and \eqref{comparison 2} in the $m=2$ case are still valid, namely 
\begin{equation}\label{comparison 1'}
\log\Big(\frac{t+1}{2}\Big) \leq G(t) \leq \log t\qquad \text{ for all } t>1
\end{equation}
and
\begin{equation}\label{comparison 2'}
e^s \leq G^{-1}(s) \leq 2e^s-1 \qquad \text{ for all } s>0.
\end{equation}
Applying to \eqref{comparison 2'} the last property in \eqref{sigma, md}, we may infer:
\begin{equation}\label{comparison 3'}
e^{\alpha s} \leq \sigma \circ G^{-1}(s) \leq 2^\alpha e^{\alpha s}\qquad \text{ for all } s>0.
\end{equation}
In addition, note that \eqref{h def} still holds true. In fact, there exists a universal constant $\delta>0$ such that for all $k \geq 1$, we can find $h=h(k)> k$ satisfying
\begin{equation}\label{h'}
[h,h+1] \subset \big[G^{-1}(k+\delta),\, G^{-1}(k+1-\delta)\big].
\end{equation}
For example, $\delta := 4^{-1}(1-\log 2)$ ensures that the length of the interval on the right-hand side of \eqref{h'} is greater than $2$.

Let the choices for $\phi_k, \alpha_k, \beta_k$ and $u_k$ be the same as in the $m=2$ case. It follows from \eqref{laplacian, Lp} and \eqref{uk, Lp} that
\begin{eqnarray}
&& \|u_k\|_{L^p}^p \leq 4(\gamma_m)^p e^{2(k+1)},\label{uk, Lp upper bd, mD}\\
&&\|\Delta u_k\|_{L^p}^p \leq (\gamma_m)^p \|\phi''\|^p_{L^\infty([0,1])} e^{-2(p-1)(k+1)},\label{laplacian uk, Lp upper bd, mD}
\end{eqnarray}
which are similar to \eqref{uk Lp upper bd} and \eqref{laplacian Lp upper bd} for $m=2$.

Now we shall specify the warping function. Again, the idea is to introduce high-frequency oscillations to $\sigma$. In view of the final line in \eqref{sigma, md}, the graph of $\sigma|[1,\infty[$ lies in the strip 
\begin{equation*}
S:= \Big\{ (t,y) \in \R^2: t \geq 1, \, t^\alpha \leq y \leq (t+1)^\alpha \Big\}.
\end{equation*}
Let us denote by 
\begin{equation*}
S_k := S \cap \big\{ k\leq t \leq k+1\big\}\qquad \text{ for each } k \in\mathbb{N}.
\end{equation*}
Note that the height of the window $S_k$ shrinks to $0$ as $k \nearrow \infty$. We introduce the parameter:
\begin{equation}
\eta_k := \min_{t \in [k,k+1]} \frac{(t+1)^\alpha - t^\alpha}{10}.
\end{equation}
As discussed above, $\eta_k \searrow 0$ as $k \nearrow \infty$. Moreover, it is easy to see that one can place a cube $Q'_k$, whose sides are parallel to the Cartesian axes and have length $\eta_k$, inside the window $S_k$. 

For $k \in \mathbb{N}$ fixed momentarily, let us define $\sigma$ on part of $[k,k+1]$. More precisely, we shall require that the graph of $\sigma$ over the horizontal projection of the cube $Q_k'$ is contained in $Q_k'$. For this purpose, we can  carry out a construction slightly simpler that in \cite{gp} for the $m=2$ case. 

Indeed, let $\sigma([z_k, z_k+\eta_k])$ be the juxtaposition of $\ell_k$ sawtooth functions of step length $$\delta_k := \frac{\eta_k}{2\ell_k}.$$ Each sawtooth function (modulo an obvious translation) increases from $0$ to $\eta_k$ in step-length $\delta_k$, and then decreases from $\eta_k$ to $0$ in another step-length $\delta_k$. Finally, we smooth out the corners by modifying on $(2\ell_k)$ intervals of the length $\delta_k^{10}$. In this way we complete the definition of $\sigma$ inside $Q_k'$; it is smooth and has gradient $|\sigma'| = \delta_k^{-1} = \nicefrac{2\ell_k}{\eta_k}$ for a large portion of the domain, {\it i.e.,} the horizontal projection of $Q'_k$. We shall specify the small parameter $\delta_k$ and the large parameter $\ell_k$ later in the proof. In passing let us note that, roughly speaking, the parameters $(\ell_k, \delta_k, \eta_k)$ play the role of $(n_k, \e_k, 1)$ as in $\S \ref{sec: m=2}$. 

In the above paragraph we defined $\sigma$ in $Q_k'$. Now let us extend it globally. For this purpose, consider a sequence $\{k_j\}_{j=1}^\infty$ which tends to $\infty$ as $j \nearrow \infty$. Let $h_j=h(k_j)$ be defined as in \eqref{h'}. As the Green's function $G$ is monotonically increasing, in view of \eqref{h'} one can choose $\{k_j\}$ so that the cubes $Q_{h_j}'$ are disjoint. Let $\sigma$ be defined in each $Q_{h_j}'$ as above. Outside these cubes we take $\sigma$ to be any smooth function satisfying the properties in \eqref{sigma, md}, and by a simple glueing argument we can obtain $\sigma \in C^\infty([0,\infty[)$. For notational convenience, in the sequel let us  relabel $k=k_j$ and $Q'_k = Q'_{h_j} \equiv Q'_{h(k_j)}$.

It remains to bound $\|\na\na u_k\|_{L^p}$ from below. First of all, by \eqref{hessian, Lp}, the choice of $\phi_k$ and the upper bound in \eqref{comparison 3'}, we have 
\begin{align*}
\|\na \na u_k\|_{L^p}^p \geq (\gamma_m)^p 2^{-\alpha p} e^{-\alpha p (k+1)}  \int_{k+\delta}^{k+1-\delta} \big|\sigma' \circ G^{-1}(s)\big|^p \,\big|\sigma \circ G^{-1}(s)\big|^{(1-m)(p-2)}\,\dd s. 
\end{align*}
Utilising once more \eqref{comparison 3'}, we get
\begin{align*}
\|\na \na u_k\|_{L^p}^p \geq (\gamma_m)^p 2^{-\alpha (p-1)} e^{-(k+1)(\alpha p + p -2)}  \int_{k+\delta}^{k+1-\delta} \big|\sigma' \circ G^{-1}(s)\big|^p\,\dd s. 
\end{align*}
For $\dim\,\M=m$ there holds $G'(r) = \sigma^{1-m}(r)$, so
\begin{equation*}
(G^{-1})'(s) = \frac{1}{G'[G^{-1}(s)]} = \sigma^{m-1}[G^{-1}(s)]. 
\end{equation*}
It yields that
\begin{align*}
\|\na \na u_k\|_{L^p}^p \geq (\gamma_m)^p 2^{-\alpha (p-1)} e^{-(k+1)(\alpha p + p -2)}  \int_{k+\delta}^{k+1-\delta} \big|\sigma' \circ G^{-1}(s)\big|^p \big|\sigma^{1-m}\circ G^{-1}(s)\big| (G^{-1})'(s)\,\dd s.
\end{align*}
Thus, changing the variables $s \mapsto r=G^{-1}(s)$ and invoking \eqref{h'}, we arrive at
\begin{align*}
\|\na \na u_k\|_{L^p}^p \geq  (\gamma_m)^p 2^{-\alpha (p-1)} e^{-(k+1)(\alpha p + p -2)} \int_h^{h+1} |\sigma'(r)|^p \sigma^{1-m}(r) \,\dd r.
\end{align*}
By \eqref{comparison 3'}, one further gets
\begin{align*}
\|\na \na u_k\|_{L^p}^p \geq (\gamma_m)^p 2^{-\alpha (p-1)} e^{-(k+1)(\alpha p + p -2)-k-\delta}\int_h^{h+1} |\sigma'(r)|^p \,\dd r,
\end{align*}
where $h=h(k)>k$ is chosen as in \eqref{h'}.

To continue, it is crucial to note that in some subinterval of $[h,h+1]$ of length $\eta_k$, $\sigma$ is highly oscillatory. This is due to our choice of $Q'_k$ and the definition of $\sigma$ thereon. More precisely, we can deduce the bound
\begin{align}\label{penultimate}
\|\na \na u_k\|_{L^p}^p &\geq (\gamma_m)^p 2^{-\alpha (p-1)} e^{-(k+1)(\alpha p + p -2)-k-\delta} (\eta_k - \delta_k^9) (\delta_k)^{-p},
\end{align}
where $\delta>0$ is universal as before. Here, recall that $\eta_k = 2\delta_k \ell_k \searrow 0$ for $\ell_k \nearrow \infty$ to be determined. We shall select some $\delta_k$ that shrinks to $0$ much more rapidly than $\eta_k \sim (k+1)^\alpha - k^\alpha = (k+1)^{\frac{1}{m-1}} - k^{\frac{1}{m-1}}$ does. Indeed, let us require that
\begin{equation*}
\begin{cases}
\delta_k^9 \leq \frac{\eta_k}{2},\\
\delta_k \leq \Big( \frac{\eta_k}{2} e^{-e^k}\Big)^{\frac{1}{p}}.
\end{cases}
\end{equation*}
The above two conditions give us
\begin{equation}\label{double exponential}
(\eta_k - \delta_k^9) (\delta_k)^{-p} \geq e^{e^{k}};
\end{equation}
while the other term on the right-hand side of \eqref{penultimate} is
\begin{equation*}
(\gamma_m)^p 2^{-\alpha (p-1)} e^{-(k+1)(\alpha p + p -2)-k-\delta} = C_3 e^{- C_4k},
\end{equation*}
with $C_3$, $C_4$ being positive constants depending only on $m$ and $p$, and with $\delta$ being a fixed universal constant as before.

To conclude the proof, we can deduce from \eqref{double exponential} and \eqref{penultimate} that for any sufficiently large $k \in \mathbb{N}$, there holds
\begin{equation*}
\|\na \na u_k\|^p_{L^p} \gtrsim  e^{k^{1000}} 
\end{equation*}
with the constants involved in $\gtrsim$ depending on $m$ and $p$. On the other hand, in \eqref{uk, Lp upper bd, mD}\eqref{laplacian uk, Lp upper bd, mD} we have already proved that
\begin{equation*}
\|u_k\|^p_{L^p} \lesssim e^k,\qquad \|\Delta u_k\|^p_{L^p} \lesssim e^{-2(p-1)k};
\end{equation*}
the constants in $\lesssim$ depending on $m$, $p$ and the $C^2$-norm of $\phi$. Finally, the choice of a large $k$ gives us the contradiction to \eqref{cz}, hence the proof of Theorem \ref{thm} for $m \geq 3$ is complete.

\section{Concluding Remark}
The open manifold $(\M,g)$ constructed in this note has no bound on the norm of the Riemann curvature, and its injectivity radius degenerates. On the other hand, if the Ricci curvature of $(\M,g)$ is bounded from the above and below, and if the injectivity radius is strictly positive, then the $L^p$-Calder\'{o}n--Zygmund estimate is valid on $(\M,g)$ for any $1<p<\infty$ ({\it cf.} Theorem $C$, \cite{gp}). It is interesting to seek for the minimal geometric boundedness assumptions on $(\M,g)$ that ensures the validity of the  $L^p$-Calder\'{o}n--Zygmund estimate.

\bigskip
\noindent
{\bf Acknowledgement}.
This work has been done during Siran Li's stay as a CRM--ISM postdoctoral fellow at Centre de Recherches Math\'{e}matiques, Universit\'{e} de Montr\'{e}al and Institut des Sciences Math\'{e}matiques. The author would like to thank these institutions for their hospitality. Siran Li also thanks Jianchun Chu for insightful discussions on problems in global analysis.

\end{document}